\documentclass{article}
\usepackage{fullpage}
\usepackage{hyperref}
\usepackage{graphicx}

\newtheorem{theorem}{Theorem}
\newtheorem{claim}{Claim}

\newtheorem{definition}{Definition}

\title{What is a Theorem?}
\author{Jeffrey C.\ Jackson \\ \\
{\small \textit{ Mathematics and Computer Science Department}}\\
{\small \textit{ Duquesne University}}\\
{\small \texttt{jacksonj@duq.edu}}}
\date{\today \\
$\Rightarrow$ \textbf{discussion of this paper is welcome at \url{drjeffjackson.blogspot.com}} $\Leftarrow$}

\begin{document}

\maketitle

\begin{abstract}
  General acceptance of a mathematical proposition $P$ as a theorem
  requires convincing evidence that a proof of $P$ exists.  But what
  constitutes ``convincing evidence?''  I will argue that, given the
  types of evidence that are currently accepted as convincing, it is
  inconsistent to deny similar acceptance to the evidence provided for
  the existence of proofs by certain randomized computations.
\end{abstract}

\section{Empirical Theorems}

\begin{figure}
\centering
\includegraphics[scale=0.30]{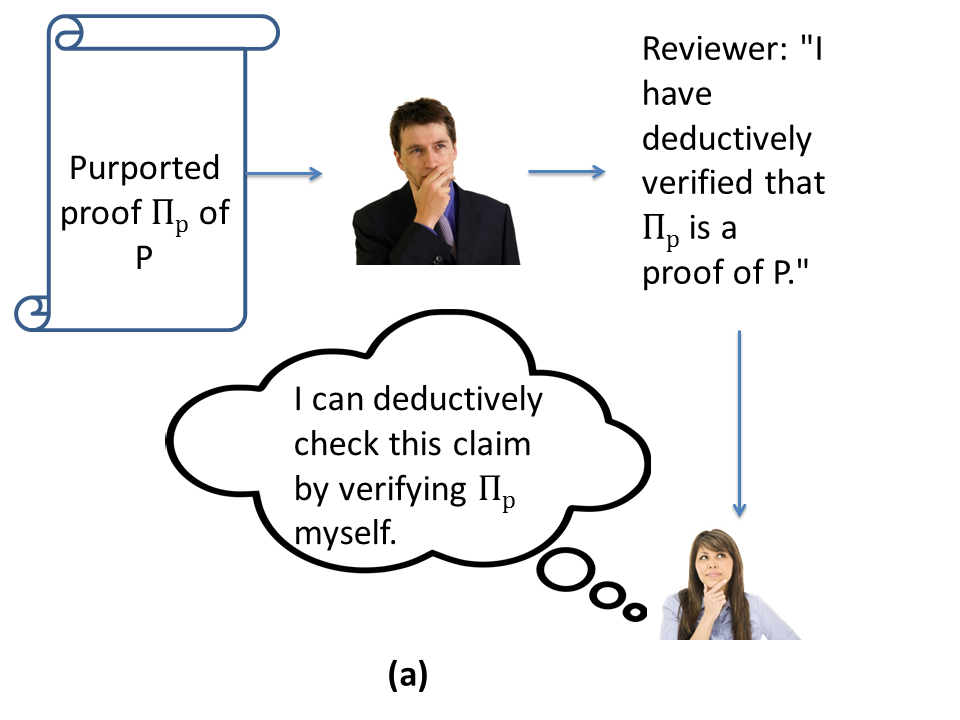}
\includegraphics[scale=0.30]{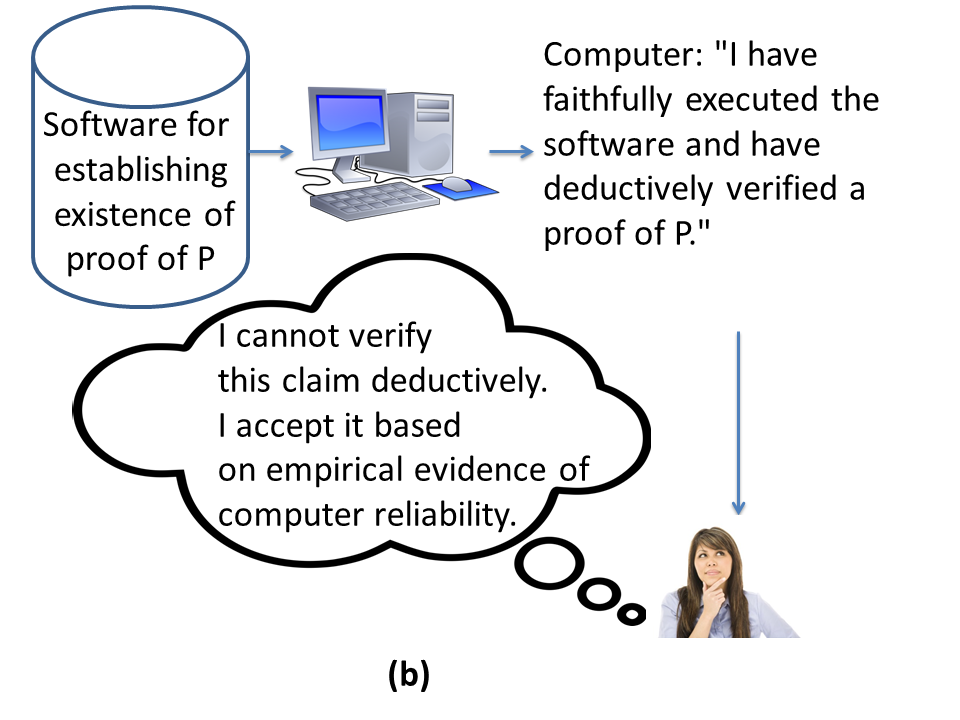}
\caption{Two approaches to accepting that the proof of a
  proposition exists. (a) Traditional. (b) Empirical.}
 \label{figaccept}%
\end{figure}

\begin{definition} \label{citable}%
  A proposition $P$ is citable as a theorem if and only if the
  mathematical community accepts that a proof $\Pi_P$ of proposition
  $P$ exists.
\end{definition}
Given that this is an acceptable definition, the question we will
focus on is, when should the mathematical community accept that a
proof of a proposition exists?  Figure~\ref{figaccept}(a) portrays the
traditional answer: A proposition $P$ is citable as a theorem if a
purported proof $\Pi_P$ of $P$ has been carefully checked through the
peer review process.  Individual mathematicians are then free to cite
$P$ as a theorem without fear of having their decision questioned.
This approach has the benefit that, at least in principle, any
mathematician (or team of mathematicians) can at any time verify for
themselves that $\Pi_P$ is indeed a proof of $P$.  This adds to
our confidence in theorems accepted in this manner, since there is
reason to believe that propositions that erroneously slip through the
peer review process will, if they are of sufficient interest to the
wider community, ultimately be discovered and retracted as theorems
(at least until a corrected proof can be offered).

In recent years, however, it has been accepted that certain proofs
exist even though no human mathematician has fully verified their
existence and no human or team of humans is likely to ever do so.
Instead, there are now propositions $P$ for which execution of
computer software has been taken as compelling evidence for the
existence of a proof $\Pi_P$.
%% The computational evidence can take
%% various forms. Sometimes, all of the reasoning behind $\Pi_P$ is
%% supplied by humans, and the computer simply carries out a massive
%% computation in order to verify that some property holds.  Other times,
%% the computer might also produce at least some of the logical steps of
%% a formal proof $\Pi_P$.  In any event, our interest is in propositions
%% $P$ for which the belief that a proof $\Pi_P$ exists necessarily
%% relies on evidence produced by computers, evidence that is thought to
%% be beyond human ability to verify.
This approach to accepting the existence of a proof of a proposition
is illustrated in Figure~\ref{figaccept}(b).

%% that are accepted as theorems
%% because a computer has produced an output affirming that a proof
%% $\Pi_P$ of $P$ exists and it is accepted that this evidence is
%% reliable.

%% has 
%% purportedly carried out computations that will result in
%% an affirmative output if and only if a proof of $P$ exists, and the
%% computer has produced that output.\footnote{It does not matter for our
%%   purposes whether the purported proof of $P$ was produced by a
%%   computer search or by traditional means.  What we are interested in
%%   is theorems $P$ for which the only evidence that a proof $\Pi_P$ of
%%   $P$ exists is evidence produced by running a computer.}  For each
%% such $P$, it is accepted by the mathematical community that a computer
%% has produced and formally verified a proof of $P$, although it is also
%% acknowledged that no human or humans will likely ever carry out their
  %% own verification.
  
Perhaps the best-known example of such a proposition $P$ is what was
long known as the Four Color Conjecture, but which is now universally
regarded as the Four Color Theorem.\footnote{There were questions for
  many years following the initial claim by Appel and Haken of a
  computer-assisted proof \cite{AppHak76} of the Four Color Theorem.
  However, there seem to be no serious qualms with the
  approach---advocated as a general technique by Hales \cite{Hal08}
  and implemented specifically for the Four Color Theorem by Gonthier
  \cite{Gon08}---of using general-purpose theorem-proving software to
  generate and check a formal proof of the theorem.}  As another
example, one that will be of some interest for purposes of this paper,
propositions of the form ``$n$ is a prime'' for very large $n$ have
been accepted based on a combination of a traditional theorem
establishing a test for Mersenne primes and massive computations
applying that test, computations that will almost certainly never be
replicated by humans; see, \textit{e.g.},
\cite{WolKur09}.  %% Although there has
%% certainly been controversy over
%% the past $40$ years or so over whether or not various computer-based
%% claims could unreservedly be considered theorems, my purpose here is
%% not to review this history but to simply note that there are some
%% propositions $P$, such as those cited, that are now beyond that
%% controversy.  
It is accepted that the Four Color Theorem can be cited as a theorem,
and it is accepted that $2^{43112609}-1$ is a prime, even though the
evidence we have that proofs of these theorems exist relies on
computational experiments.

I will refer to propositions that have been accepted as theorems based
in whole or in part on the approach of Figure~\ref{figaccept}(b) as
\emph{empirical theorems}.  For as the figure illustrates, we cannot
accept an empirical theorem without first concluding that a computer
has verified the existence of a proof.  And what justifies this
conclusion?  Without a doubt, deductive reasoning can play a role: We
might have used formal methods to develop the software and design the
hardware, there might be traditional theorems underlying the
algorithms that are employed, and so on.  However, in the end, even if
we have good reason to be convinced that the software is perfectly
written and the hardware design is flawless, a computer program is
executed by a physical device that might behave differently than we
would expect for a variety of reasons (material defects, tampering, an
electromagnetic pulse, etc.).  How then can we ``know'' that an
execution of a computer program has not produced misleading output?
Ultimately, it is scientific, not mathematical, induction that leads
us to conclude, on the basis of past experience, that modern digital
computers are extremely reliable and that it is therefore reasonable
to trust that they have faithfully executed their instructions.  This
implies, as Figure~\ref{figaccept}(b) illustrates, that if we question
the fidelity of any given execution of a theorem-verifying system, our
only viable recourse is to run another physical experiment producing
further empirical evidence.

\begin{figure}
\centering
\includegraphics[scale=0.33]{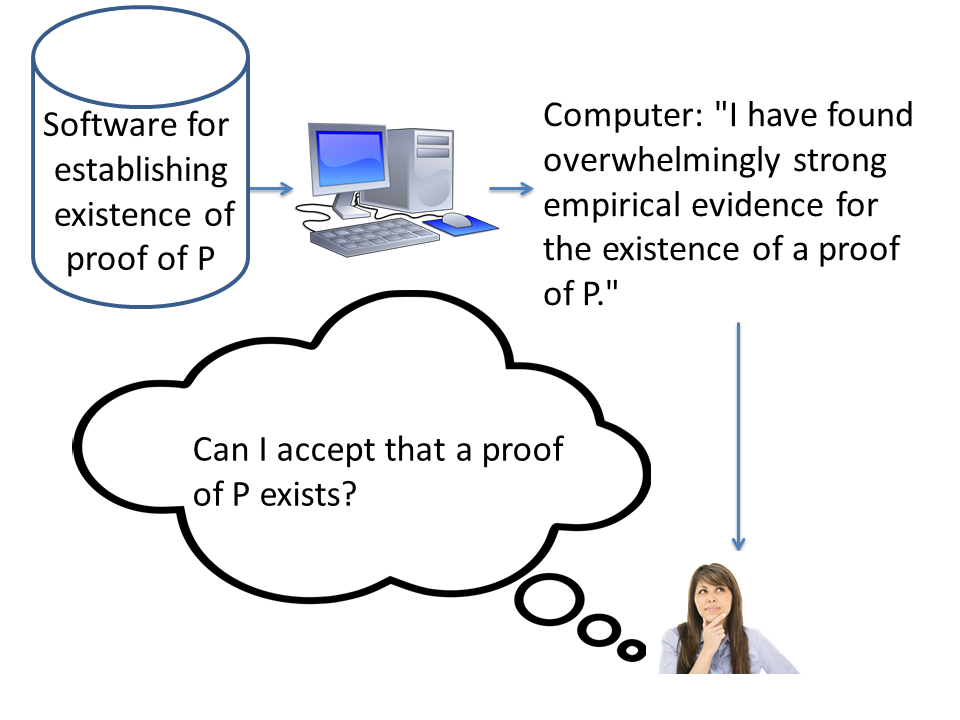}
\caption{A third approach to accepting that the proof of a
  proposition exists.}
 \label{two}%
\end{figure} 

The question considered here is whether we should recognize a third
method for accepting the existence of proofs (Figure~\ref{two}).  In
this method, rather than requiring the computer to produce and
deductively verify a proof of $P$, the computer is allowed to employ
empirical reasoning to reach the conclusion that a proof $P$ exists.
%We will begin exploring this question with a case study.

\section{A Case Study: Primality Testing}

Primality testing is one area that provides examples of empirical
theorems.  Published primality theorems, such as \cite{WolKur09}, are
based on \textit{deterministic} algorithms for primality testing.  On
a fixed input, if such an algorithm is faithfully executed repeatedly
then it will produce the same output on every execution.

However, deterministic algorithms are not the only approach to
primality testing.  For instance, the Miller-Rabin primality testing
algorithm \cite{Rab80} has been widely used in cryptographic
applications and elsewhere.  This algorithm is \textit{randomized} in
the sense that the algorithm is provided not only the input that a
deterministic algorithm would receive---an integer to be tested for
primality, in the case of Miller-Rabin---but also with a stream of
bits thought to be random.\footnote{The bits typically treated as
  ``random'' in computer software are actually pseudo-random,
  generated by numerical processes that, although deterministic,
  produce bit sequences that pass various statistical tests for
  randomness.  However, it is possible to obtain bit streams produced
  by physical processes that are generally accepted as behaving
  randomly.  For instance, at the time of this writing, bits from a
  quantum-mechanical process are being generated and made available
  online by the Australian National University
  (\url{https://qrng.anu.edu.au/}).}  The output of the algorithm
depends on both the input and the random bits provided, so repeated
faithful executions of the algorithm on a fixed input but with
different bit streams can produce differing outputs.  The question we
will focus on in this section is, what sort of evidence for the
existence of a proof of primality does execution of Miller-Rabin
provide?  In the next section, we will compare Miller-Rabin evidence
with the sort of evidence provided by a deterministic primality
testing algorithm.

Let us begin by defining a key component of the Miller-Rabin
algorithm, the \textit{witness to compositeness} function $W_n(b)$ for
$1 \leq b \leq n$. $W_n(b)$ outputs ``composite'' if $b^{n-1} \not
\equiv 1\pmod{n}$ or there exists an $i$ such that $2^i | (n-1)$ and
$\gcd(b^{(n-1)/2^i}-1,n)$ is not $1$ or $n$.  Otherwise, $W_n(b)$
outputs ``indeterminate.''  This function would clearly be simple to
program, and in fact it is similar in form to the Lucas-Lehmer test
used to verify that $2^{43112609}-1$ is prime \cite{WolKur09}.
$W_n(b)$ has the following properties:
\begin{theorem}[Miller]
If $n$ is prime then $W_n(b)$ will return ``indeterminate'' for
every $1 \leq b < n$.
\end{theorem}
\begin{theorem}[Rabin]
If $n>4$ is composite then $W_n(b)$ will return ``composite'' for
at least $\frac{3}{4}$ of the values $1 \leq b < n$.
\end{theorem}

These properties of $W_n(b)$ imply that the following
deterministic approach could, in principle, produce a proof of the
primality of $n$: Select $\lfloor \frac{(n-1)}{4} \rfloor + 1$
distinct integer values in $[1,n)$ and execute $W_n$ on each value
until either the output is ``composite'' or $W_n$ has been applied to
every value.  If the output is ever ``composite,'' $n$ is composite.
Otherwise, $n$ is prime.

This approach, although deterministic, would be utterly impractical
for large $n$.  Thus, instead, the Miller-Rabin algorithm selects some
number $k$ of integers in $[1,n)$ uniformly at random (using its
  stream of random bits) and executes $W_n$ on each of these integers.
  If any output of $W_n$ on one of the $k$ integers is ``composite,''
  the algorithm output is ``composite,'' and otherwise the algorithm
  output is ``prime.''  It is easy to see that if the algorithm is
  faithfully executed, on prime $n$ it will definitely output
  ``prime'' and on composite $n$ it will output ``composite'' with
  probability, over the random choice of the $k$ test integers, at
  least $1-1/4^k$.

Although it might seem at first that when Miller-Rabin outputs
``prime'' it is merely providing evidence for the primality of a
number $n$, it is important to recognize that it is providing more
than this: It is providing evidence that a proof of the primality of
$n$ exists.  Specifically, it is providing evidence that, were we to
continue executing $W_n$ faithfully on sufficiently many integer
values, we would produce a deterministic proof of $n$'s primality.
This distinction between evidence for a proposition $P$ and evidence
for the existence of a proof $\Pi_P$ of $P$ is important because, per
Definition~\ref{citable}, $P$ is only a citable theorem if it is
accepted that a proof $\Pi_P$ exists.  If Miller-Rabin merely provided
evidence for the primality of $n$, this evidence might increase our
confidence in a \textit{conjecture} of $n$'s primality, but we would
be no closer to establishing a citable \textit{theorem} stating that
$n$ is prime.
% But, again, Miller-Rabin actually provides evidence for
% the existence of a proof of primality, which in turn provides evidence
% that ``$n$ is prime'' is a theorem.

% Finally, let us consider the sort of evidence we would have were we to
% physically execute Miller-Rabin for a given $n$ with a reasonably
% large value for $k$, say, $50$.  If the execution were to output
% ``prime,'' this would provide evidence of the existence of a proof of
% the primality of $n$.  Admittedly, there would be a (slim) chance that
% this output would be faulty due to providing the algorithm with an
% ``unlucky'' stream of random bits.  And there would also be a (slim)
% chance that this output would be faulty due to a malfunction of some
% aspect of the physical system executing the algorithm.  Note, however,
% that comparing this description of the evidence produced by an
% execution of Miller-Rabin with Assertion~\ref{assertion:one}, we see
% that theorems have been accepted based on comparable evidence of the
% existence of proof.  We pursue this line of reasoning in more detail
% in the next section.

\section{Comparing Deterministic and Randomized Evidence of Proof}

Imagine now that Miller-Rabin has been executed with large integer
input $n$ on a digital computer with a number of tests $k$ large
enough that the probability of error due to randomization is
comparable with the probability that an execution of a deterministic
primality testing algorithm on a digital computer will mistakenly
output ``prime'' when its input is actually composite.  Thus, as noted
in the previous section, this execution of Miller-Rabin produces
extremely strong empirical evidence for the existence of a proof of
the proposition ``$n$ is prime.''  Should ``$n$ is prime'' be citable
as a theorem?

Currently, the consensus answer of the mathematical community seems to
be a clear, even emphatic, ``no.''  I will argue that consistency
suggests that the answer should instead be, ``yes.''  The argument
will be familiar to computer scientists, as something like it has long
been used to justify the use of randomized algorithms in computing
practice.  I hope to frame the argument in a way that makes it clear
to mathematicians that randomized algorithms can in some cases be used
to establish the existence of (deductive) proofs of theorems.

\begin{figure}
\includegraphics[scale=0.33]{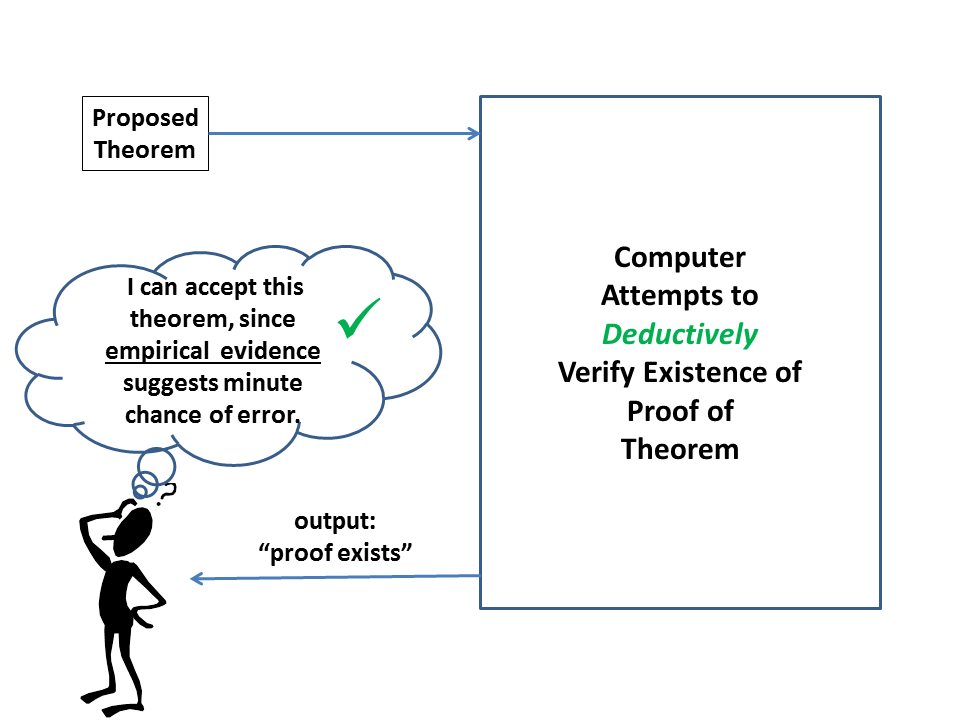}
\includegraphics[scale=0.33]{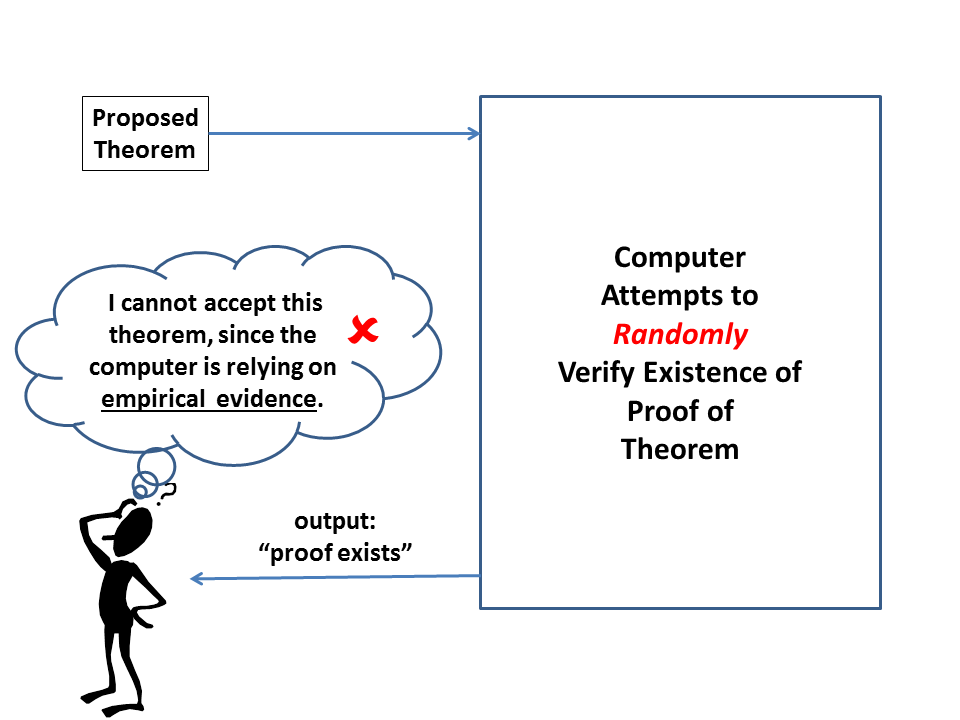}
\caption{Mathematician Mel considering the results of two different
  computational approaches to verifying the existence of a proof of a
  theorem.}
 \label{figcompare}%
\end{figure}

The argument, in a nutshell, is this: Let us imagine a mathematician
Mel who states, ``I can accept an empirical theorem, but only if I
believe that a computer has produced and deductively checked a
proof of the theorem.''  The situation is illustrated in
Figure~\ref{figcompare}.  If we view the computer as a sort
of assistant in the process of establishing a theorem, then
the consistency problem is this: 
% Mel uses
% empirical evidence to arrive at his belief that the computer has
% deductively checked a proof.  That is, 
Mel is allowed to ignore the very small chance of being misled by
empirical evidence but begrudges the assistant exercising the same
privilege.  Put another way, although scientific inference plays a
fundamental role in Mel's accepting an empirical theorem, Mel wants to
deny that such inference is suitable in a secondary role.

To be sure, it is reasonable that Mel might \textit{prefer} having an
empirical theorem based on a deterministic algorithm to having one
based on a randomized algorithm, much as a traditional theorem might
be preferred to an empirical theorem or a constructive proof to an
existence proof.  But if empirical evidence for the existence of a
proof is going to be accepted---and it is---then we are employing
something like the scientific method, which would seem to imply that
it should be the confidence we have in the evidence produced, and not
the form of the experiments used to produce the evidence, that
determines whether or not we accept the evidence.

I hope that at this point the reader will allow that, at least at a
first glance, it appears that it indeed might be difficult to argue
that it is consistent for Mel---and the mathematical community in
general---to accept empirical theorems but to reject out of hand
randomized empirical evidence for the existence of proofs of such
theorems.  This leads me to make the following claim, which I hope
will lead to further discussion:
\begin{claim}
  There is some reasonably small integer $k$ such that, if $n$ is
  fixed and the Miller-Rabin witness to compositeness function $W_n$
  is correctly implemented and executed by a presumed-reliable
  computer on $k$ uniform-randomly selected test integers and outputs
  ``indeterminate'' on every test, then it should be accepted that a
  proof of the theorem ``$n$ is prime'' exists.  More generally, if
  execution of any randomized algorithm provides a comparable level of
  evidence for the existence of a proof of a proposition $P$, then it
  should be accepted that a proof of $P$ exists.
\end{claim}
\noindent
It is beyond the scope of this paper to argue for any particular value
for $k$.

\section{Summary and Future Work}

I have noted that empirical theorems are currently accepted by the
mathematical community and have emphasized that non-deductive,
empirical reasoning is essential in accepting these theorems.  Based
on these observations, I have claimed that it is inconsistent to
disallow empirical elements in algorithms designed to establish the
existence of proofs.  My hope is that this argument will spur
discussion that will, sooner rather than later, lead the mathematical
community to embrace the notion that certain randomized algorithms can
legitimately play a role in establishing theorems.

If this notion were to be embraced in concept, before putting the
randomized approach into practice it would also be necessary to agree
on an acceptable error threshold for randomized algorithms used in
support of establishing theorems.  A possible starting point would be
to attempt to quantify the chance that computer execution of a
deterministic algorithm will produce a misleading result.

\section*{Acknowledgments}

Brad Lucier and Joel Hass each provided a number of helpful comments
on earlier drafts of this paper, and Brendon LaBuz, Alex Lipecky, and
Rachael Neilan asked questions that stimulated my thinking.

\bibliographystyle{plain}
\bibliography{../references}

\begin{thebibliography}{1}

\bibitem{AppHak76}
K.~Appel and W.~Haken.
\newblock Every planar map is four colorable.
\newblock {\em Bull. Amer. Math. Soc.}, 82(5):711--712, 09 1976.

\bibitem{Gon08}
Georges Gonthier.
\newblock Formal proof---the four-color theorem.
\newblock {\em Notices of the AMS}, 55(11):1382--1393, December 2008.

\bibitem{Hal08}
Thomas~C Hales.
\newblock Formal proof.
\newblock {\em Notices of the AMS}, 55(11):1370--1380, 2008.

\bibitem{Rab80}
Michael~O Rabin.
\newblock Probabilistic algorithm for testing primality.
\newblock {\em Journal of Number Theory}, 12(1):128 -- 138, 1980.

\bibitem{WolKur09}
George Woltman and Scott Kurowski.
\newblock On the discovery of the 45th and 46th known mersenne primes.
\newblock {\em Fibonacci Quarterly}, 46/47(3):194--197, August 2009.

\end{thebibliography}

\end{document}